\documentstyle[12pt]{article}
\newcommand{\fr}{ \frac}
\newcommand{\beq}{\begin{equation}}
\newcommand{\eeq}{\end{equation}}
\newcommand{\br}{\begin{array}}
\newcommand{\er}{\end{array}}
\newcommand{\lb}{\label}
\newcommand{\r}{\ref}
\newcommand{\ar}{\rightarrow}
\newcommand{\dl}{\delta}
\newcommand{\ov}{\overline}
\newcommand{\lan}{\langle}
\newcommand{\rn}{\rangle}
\newcommand{\pr}{\prime}

\newcommand{\lm}{\lambda}

\newcommand{\var}{\epsilon}

\begin{document}
\begin{center}
{\large \bf Representations of $SU(1,1)$ in Non-commutative Space 
Generated by the Heisenberg Algebra}
\end{center}

\vspace{1cm}

\begin{flushleft}
H. Ahmedov$^1$  and I. H. Duru$^{2,1}$
\vspace{.5cm}

{\small 
1. Feza G\"ursey Institute,  P.O. Box 6, 81220,  \c{C}engelk\"{o}y, 
Istanbul, Turkey 
\footnote{E--mail : hagi@gursey.gov.tr and duru@gursey.gov.tr}.

2. Trakya University, Mathematics Department, P.O. Box 126, 
Edirne, Turkey.}
\end{flushleft}
\vspace{1cm} \noindent 
{\bf Abstract}: $SU(1,1)$ is considered as the automorphism group 
of the Heisenberg algebra $H$. 
The basis in the Hilbert space $K$ of functions on $H$ on which 
the irreducible representations of the group are realized 
is explicitly constructed. The addition theorems  are  derived. 

\begin{center}
March 2000
\end{center}

\noindent
{\bf 1. Introduction}

\vspace{5mm}
\noindent
Investigating the properties of manifolds by means of the symmetries they 
admit has a long history. Non-commutative geometries have become the subject 
of similar studies in recent decades. For example there exists an extensive 
literature on the $q$-deformed groups $E_q(2)$ and $SU_q(2)$ which are the 
automorphism groups of the quantum plane $zz^* =qz^*z$ and the quantum sphere 
respectively \cite{w}. Using group theoretical methods the invariant distance 
and the Green functions have also been written in these deformed spaces
\cite{a}.

In the recent work we started to analyze yet another non-commutative space 
$[z, z^*]=1$ ( i. e. the space generated by the Heisenberg  algebra )
by means of its automorphism groups: We  considered $E(2)$ group 
transformations in $z, z^*$ space; and constructed the basis 
(which are written in terms of the Kummer functions) in this space 
where the  unitary irreducible representations of $E(2)$ are realized 
\cite{a1}. This analysis revealed a peculiar connection between the 
$2$-dimensional Euclidean group and the Kummer functions.

In the present work we continue to study the same non-commutative space 
$[z, z^*]=1$, this time by means of the other admissible automorphism group 
$SU(1,1)$.

In Section 2 we define $SU(1,1)$ in the Heisenberg algebra $H$ and construct 
the unitary representations of the group in the Hilbert space $X$ where 
$H$ is realized.

In Section 3 we classify the invariant subspaces in the space of the bounded 
functions on $H$ where the  irreducible representations of $SU(1,1)$ are 
realized. 

In Section 4 we show that in the Hilbert space $K$ of the square integrable 
functions  only principal series is unitary. We construct the orthonormal 
basis in $K$ which can  be written in terms of the Jacobi functions.

Section 5 is devoted to the addition theorems.
These theorems provide a group theoretical interpretation for the already 
existing identities involving the hypergeometric functions  which all are 
actually the Jacobi functions. They may also lead to new identities. 

\vspace{2cm}

\noindent
{\bf 2. Weyl representations of $SU(1,1)$ }

\vspace{5mm}
\noindent
The one dimensional Heisenberg algebra $H$ is the 
3-dimensional vector space with the basis elements $\{ z, z^*, 1\}$ 
and the bilinear antisymmetric product  
\beq\lb{com}
[z, z^*]=1.
\eeq
The $*$-representation of $H$ in the suitable dense subspace of the Hilbert 
space $X$ with the complete orthonormal basis $\{\mid n\rn \}$, 
$n=0,1,2, ...$ is given by 
\beq\lb{r1}
z \mid n\rn =\sqrt{n} \mid n-1\rn, \ \ \ 
z^*\mid n\rn =\sqrt{n+1} \mid n+1\rn. 
\eeq
Let us represent the pseudo-unitary group $SU(1,1)$ in the vector space $H$ 
\beq\lb{tr}
g\left(
\begin{array}{c} 
z \\ 
z^* 
\end{array} 
\right) =
\left(
\begin{array}{cc} 
a & b  \\
\ov{b} & \ov{a}  
\end{array} 
\right) 
\left(
\begin{array}{c} 
z \\ 
z^* 
\end{array} 
\right). 
\eeq
Due to 
\beq
a\ov{a}-b\ov{b}=1
\eeq
the transformations (\r{tr}) preserve the commutation relation
\beq
[g z, gz^*] =[z, z^*].
\eeq
Therefore 
\beq\lb{u}
g z = U(g) z U^{-1}(g), \ \ \ g z^* = U(g) z^* U^{-1}(g)
\eeq
where $U(g)$ is the unitary representation of $SU(1,1)$ in  $X$: 
\beq
U(g_1)U(g_2)=U(g_1g_2), \ \ \ U^*(g)=U^{-1}(g) = U(g^{-1}).
\eeq 
The Cartan decomposition for the group reads
\beq\lb{car}
g=k(\phi)h(\alpha )k(\psi), 
\eeq
where
\beq
k (\psi)=\left(
\begin{array}{cc} 
e^{\fr{i\psi}{2}} & 0  \\
0 & e^{-i\fr{\psi}{2}} 
\end{array} 
\right), \ \ \  \
h (\alpha)=\left(
\begin{array}{cc} 
\cosh \fr{\alpha}{2} & \sinh \fr{\alpha}{2}  \\
\sinh \fr{\alpha}{2} & \cosh \fr{\alpha}{2}
\end{array} 
\right).
\eeq
For the subgroup $k(\psi )$ we have
\beq
U(k(\psi ))\mid n \rn = e^{-i\fr{n\psi }{2}} \mid n\rn.
\eeq
Let us choose the following realizations for $z$, $z^*$ and $X$:
\beq
z=\fr{1}{\sqrt{2}} (x +\fr{d}{dx}), \ \ \ 
z^*=\fr{1}{\sqrt{2}} (x -\fr{d}{dx}),
\eeq
\beq
\lan x\mid n \rn = \Psi_n (x), \ \ \
\Psi_n (x) = \sqrt{\fr{e^{-x^2}}{2^n n! \sqrt{\pi} }} H_n (x),
\eeq
where $H_n$ is the Hermite polynomial. 
From 
\beq
h(\alpha ) z = \fr{1}{\sqrt{2}} (xe^{\fr{\alpha}{2}} + e^{-\fr{\alpha}{2}}\fr{d}{dx})
\eeq
and 
\beq
\int_{-\infty}^\infty dx \ov{\Psi_m(x)} \Psi_n(x) =\dl_{nm}
\eeq
we get 
\beq
U(h(\alpha)) \Psi_m(x) = e^{\fr{\alpha}{4}} \Psi_m( e^{\fr{\alpha}{2}}x).
\eeq
Matrix elements of $U(h(\alpha))$ in the basis $\{\mid n\rn \}$ reads
\beq
U_{mn}(h) \equiv \lan m\mid U(h(\alpha)) \mid n\rn = 
e^{\fr{\alpha}{4}}
\int_{-\infty}^\infty dx \ov{\Psi_m(x)} \Psi_n(e^{\fr{\alpha}{2}}x).
\eeq
Evaluating this integral  we get 
\beq
U_{m n}(h ) = \fr{2^{\fr{m-n}{2}}}{(\fr{n-m}{2})!} 
\sqrt{\fr{n!\sinh^{n-m} \fr{\alpha}{2}}
{m!\cosh^{n+m+1} \fr{\alpha}{2}}} 
F(-\fr{m}{2}, \fr{1-m}{2};1+\fr{n-m}{2}; -\sinh^2 \fr{\alpha}{2})
\eeq
if $n\geq m$ and $n+m$ is even and
\beq
U_{mn}(h) = 0
\eeq
if $n+m$ is odd. For $m\geq n$ one has to replace $m$, $n$ and $\alpha$ 
in the above formulas by $n$, $m$ and $-\alpha$ respectively.

\vspace{1cm}
\noindent
{\bf 3. Irreducible representations of $SU(1,1)$ in $H$ }

\vspace{5mm}
\noindent
The formula 
\beq\lb{rep}
T(g) F(z) = F(gz)
\eeq
defines the representation of $SU(1,1))$ in the space $K_0$ of bounded 
operators in the Hilbert space $X$ representable as the finite sums
\beq\lb{set}
F = \sum (f_n (\zeta )z^n + z^{*n}f_{-n}(\zeta)). 
\eeq
Here  $f_n(\zeta)$ are functions of the self-adjoint operator 
$\zeta =z^*z$. 
Using (\r{u}) we can rewrite (\r{rep}) in the form
\beq\lb{rep1}
T(g) F(z) = U(g)F(z) U^*(g)
\eeq
With the one parameter subgroups $g_1= h(\epsilon)$, 
$g_2= k(\fr{\pi}{2})h(\epsilon)k(-\fr{\pi}{2})$ and 
$g_3=k(\epsilon)$ of $SU(1,1)$ we associate the linear operators
$E_k: K_0\ar K_0$
\beq\lb{inf}
E_k (F) = \lim_{\epsilon\ar 0} \fr{1}{\epsilon}(T(g_k)F - F) 
\eeq
with  the  limit being taken in the strong operator topology.  
Inserting (\r{rep1}) into (\r{inf}) we get 
( with $H_\pm = -E_1 \mp iE_2, \ \ \ H = iE_3$ )
\beq\lb{real}
H_- (F) = \fr{1}{2}[F, z^2] , \ \ \ H_+ (F) = \fr{1}{2}[z^{*2}, F], \ \ \
H (F) = \fr{1}{2}[\zeta, F] ,  
\eeq
which implies the Lie algebra of $SU(1,1)$
\beq
[H_+, H_- ]= 2H, \ \ \ [H, H_\pm]= \pm H_\pm.
\eeq
The irreducible representations labelled by pair $(\tau, \var)$, 
$\tau\in C$ and $\var =0, \fr{1}{2}$ are given by the formulas \cite{v}
\beq\lb{1}
H_- D_k^{(\tau,\var)}=-(k+\tau +\var)D_{k-1}^{(\tau,\var)}, 
\eeq 
\beq\lb{2}
H_+ D_k^{(\tau,\var)}=(k-\tau +\var)D_{k+1}^{(\tau,\var)}, 
\eeq 
\beq\lb{3}
H D_k^{(\tau,\var)}=(k+\var)D_k^{(\tau,\var)}.
\eeq
(\r{real}) and (\r{3}) imply 
\beq\lb{a1}
D_k^{(\tau,\var)} = z^{*2(k+\var )}f_k^{(\tau,\var)} (\zeta) 
\eeq
for $k\geq 0$ and
\beq\lb{a2}
D_k^{(\tau,\var)}= f_k^{(\tau,\var)}  (\zeta) z^{-2(k+\var)}  
\eeq
for $k< 0$. 
By substituting (\r{a1}) in (\r{1}) and (\r{2}) with
\beq
f_k^{(\tau,\var)} (\zeta) = \sum_{n=0}^\infty \fr{(-)^n2^{n+k+\var}}{n!} 
C_{kn}z^{*n}z^n
\eeq
we get the recurrence relations
\beq
nC_{kn-1} +\fr{k+\var +\tau}{2k+2\var+n-1}C_{k-1n}-(2k+2\var +n)C_{kn}=0,
\eeq
\beq
C_{kn+1} -C_{kn+2}-(k+\var -\tau)C_{k+1n}=0,
\eeq
which are solved by
\beq
C_{kn}= \fr{\Gamma (1+\tau+\var+k+n)}{\Gamma(1+2\var+2k+n)}.
\eeq
Using
\beq
z^{*n}z^n = \zeta (\zeta -1)... (\zeta - n+1) 
\eeq
for $k\geq 0$ we get 
\beq
f_k^{(\tau,\var)} (\zeta) = (-2)^{k^\pr}
\fr{\Gamma (1+\tau +k^\pr )}
{\Gamma (1+2k^\pr)} 
F (-\zeta , 1+\tau+k^\pr; 1+2k^\pr; 2),
\eeq
where $k^\pr =k+\var$. 
The functions $f_k^{(\tau,\var)}$ for $k< 0$ is shown to be defined from the 
expression
\beq
f_k^{(\tau,\var)} (\zeta) = f_{-k}^{(\tau,-\var)} (\zeta).
\eeq
From (\r{1}), (\r{2}) and (\r{3}) we conclude 
that $SU(1,1)$ admits the following irreducible representations:

i) $T_{(\tau,\var)}: \ \  (\tau +\var) \ov{\in} Z$

ii) $T^\pm_{(\tau,\var)}: \ \  (\tau +\var) \in Z , \ 
\tau -\var < 0$, that is 
$\tau = -\fr{1}{2}, -1, -\fr{3}{2}, . . . $

iii) $T^0_{(\tau,\var)}: \ \  (\tau +\var) \in Z , \ \tau -\var \geq 0$, 
that is $\tau = 0, \fr{1}{2}, 1, \fr{3}{2}. . . $

\vspace{2mm}
\noindent
The corresponding invariant subspaces are:

i) $V_{(\tau,\var)}$  generated by 
$\{D^{(\tau,\var})_k\}_{k=-\infty}^\infty$

ii) $V^+_{(\tau,\var)}$ and  $V^-_{(\tau,\var)}$ generated by 
$\{D^{(\tau,\var})_k\}_{k=-\infty}^{\tau-\var}$ and
$\{D^{(\tau,\var})_k\}_{k=-\tau-\var}^\infty $

iii) $V^0_{(\tau,\var)}$ generated by 
$\{D^{(\tau,\var})_k\}_{k=-\tau-\var }^{\tau-\var}$

\vspace{1cm}
\noindent
{\bf 4. Unitary irreducible representations of $SU(1,1)$ in $H$ }

\vspace{5mm}
\noindent
We can define the norm in the subspace of $K_0$ with $f_n(\zeta)$ in (\r{set}) 
being  the functions with finite support in 
$Spect (\zeta ) =\{0, 1, 2, ...\}$ as
\beq
\mid \mid F\mid\mid =\sqrt{tr (F^*F)}.
\eeq
Completion of this subspace leads to the Hilbert space  $K$ of the square 
integrable functions in the linear space $H$ with the scalar product 
\beq\lb{sc}
(F,G)=tr (F^*G).
\eeq
Using (\r{rep1}), the unitarity of $U(g)$ and the property of the trace we 
conclude that the representation $T(g)$ in $K$ is unitary.
(\r{real}) implies the real structure in the Lie algebra
\beq
H_\pm^* = - H_\mp, \ \ \  H^* =  H.
\eeq
To investigate the unitarity of the irreducible representations 
in the Hilbert space $K$ classified  in the previous section   we consider 
the orthogonality condition for the basis elements $D^{(\tau, \var)}_k$. 
Using (\r{r1}) and (\r{a1}) we get 
\beq
(D^{(\tau, \var)}_k, D^{(\tau^\pr, \var^\prime)}_m ) = 
\dl_{mk}\dl_{\var\var^\pr} \sum_{n=0}^\infty \fr{(n+2k+2\var)!}{n!}
\ov{f_k^{(\tau,\var)} (-n)} f_k^{(\tau^\pr,\var)} (-n).
\eeq
Putting
\beq
s=1-e^{-t}, \ \ \ \lm =1+2(k+\var) +\mu
\eeq
in the formula \cite{e1}
\begin{eqnarray}
\sum_{n=0}^\infty \fr{\Gamma (n+\lm)}{n!\Gamma(\lm)} s^n 
F(-n, a; \lm ;2)F(-n, b; \lm ;2) =  \nonumber \\
= (1-s)^{a+b-\lm} (1+s)^{-a-b} F(a,b; \lm ; \fr{4s}{(1+s)^2})
\end{eqnarray}
and taking first the limit $\mu\ar +0$ and then  $t\ar \infty$ we 
obtain  for $\tau =-\fr{1}{2}+i\rho$, $\rho\in R$ the orthogonality relations
\beq
(D^{(-\fr{1}{2}+i\rho, \var)}_k, D^{(-\fr{1}{2}+i\rho^\pr, \var^\prime)}_m ) 
= \dl_{mk}\dl_{\var\var^\pr} \dl(\rho -\rho^\pr).
\eeq
In the deriving of the above relation we used 
\beq
F(a,b;c;1) = \fr{\Gamma (c) \Gamma (c-a-b)}{\Gamma (c-a) \Gamma (c-b)}
\eeq
and the representation  
\beq
\lim_{t\ar \infty} \fr{e^{-izt}} {z+i0} =-2\pi i\delta (z)
\eeq
for the Dirac delta function. For other values of $\tau$ there is 
no orthogonality condition. Thus in $K$ only the representation  
$T_{(\tau, \var)}$ with $\tau = -\fr{1}{2}+i\rho$ of Section 3,
which is the principal series is unitary. 

\vspace{1cm}
\noindent

\noindent
{\bf 5. The addition theorems}

\vspace{5mm}
\noindent
({\bf i}) Restriction of (\r{rep}) on the subspace $V_{(\tau, \var)}$ reads: 
\beq
T(g) D^{(\tau, \var)}_{k} = 
\sum_{n=-\infty}^\infty t_{nk}^{(\tau, \var)} (g) D^{(\tau, \var)}_{n}
\eeq
or
\beq\lb{51}
U(g) D^{(\tau, \var)}_{k} U^*(g) = 
\sum_{n=-\infty}^\infty t_{nk}^{(\tau, \var)} (g) D^{(\tau, \var)}_{n}
\eeq
where 
\begin{eqnarray}
t_{kn}^{(\tau, \var)} (g)=
\fr{e^{-i(k+\var)\phi -i(k+\var)\psi}}{(k-n)!}
\fr{\Gamma (1+ \tau - \var -n ) \sinh^{k-n} \fr{\alpha}{2}}
{\Gamma (1+\tau - \var - k ) \cosh^{k+n+2\var} \fr{\alpha}{2} } 
 \times \nonumber \\ 
\times F(-\tau-\var - n , 1+\tau-\var -n ; 1+ k-n ; -\sinh^2 \fr{\alpha}{2} )
\end{eqnarray}
are the matrix elements of the irreducible representations which are valid 
for $k \geq n$. For $k < n$ one has to replace $k$ and 
$n$ on the right hand side by $-k$ and $-n$ respectively. 

\vspace{5mm}
\noindent
({\bf ii}) Restriction of (\r{rep}) on the subspaces $V^+_{(\tau, \var)}$ and 
$V^-_{(\tau, \var)}$ gives the following addition theorems: 
\beq\lb{52}
U(g) D^{(\tau, \var)}_{k} U^*(g) = 
\sum_{n=-\infty}^{\tau -\var} t_{nk}^{(\tau, \var)} (g) D^{(\tau, \var)}_{n}
\eeq
and 
\beq\lb{53}
U(g) D^{(\tau, \var)}_{k} U^*(g) = 
\sum_{n=-\tau-\var}^\infty t_{nk}^{(\tau, \var)} (g) D^{(\tau, \var)}_{n}.
\eeq

\vspace{5mm}
\noindent
({\bf iii}) On the subspaces $V^0_{(\tau, \var)}$ the 
addition theorem reads 
\beq\lb{54}
U(g) D^{(\tau, \var)}_{k} U^*(g) = 
\sum_{n=-\tau -\var }^{\tau -\var} t_{nk}^{(\tau, \var)} (g) 
D^{(\tau, \var)}_{n}.
\eeq
Sandwiching  both sides of (\ref{51}), (\ref{52}), (\ref{53}) and (\ref{54})
between  the states $\lan l\mid $ and $\mid s\rn$ we get 
\beq\lb{5a}
\sum_{m,t=0}^\infty U_{lm}(g) \ov{U_{st}(g)} (D^{(\tau, \var)}_{k})_{mt} 
= \sum_n t_{nk}^{(\tau, \var)} (g) (D^{(\tau, \var)}_{n})_{ls}
\eeq
Multiplying (\ref{51}), (\ref{52}), (\ref{53}) and (\ref{54})  by $U(g)$ 
from the right and sandwiching them  between  the states $\lan l\mid $ 
and $\mid s\rn$ we get 
\beq\lb{5b}
\sum_{m=0}^\infty U_{lm}(g) (D^{(\tau, \var)}_{k})_{ms} 
=  \sum_{m=0}^\infty
\sum_n t_{nk}^{(\tau, \var)} (g) (D^{(\tau, \var)}_{n})_{lm}
U_{ms}(g)
\eeq
Multiplying (\ref{51}), (\ref{52}), (\ref{53}) and (\ref{54})  by 
$U^*(g)$ and $U(g)$ from the left and right respectively and sandwiching 
them  between  the states $\lan l\mid $ and $\mid s\rn$ we get 
\beq
(D^{(\tau, \var)}_{k})_{ls} 
=  \sum_{m,t=0}^\infty
\sum_n t_{kn}^{(\tau, \var)} (g) 
U_{ts}(g)\ov{U_{ml}(g)}(D^{(\tau, \var)}_{n})_{mt}
\eeq
where 
\beq
(D^{(\tau, \var)}_{k})_{mt}= 
\sqrt{\fr{m!}{t!}} f_k^{(\tau,\var}(t)\dl_{m,t+2k+2\var}
\eeq
for $k\geq 0$ and 
\beq
(D^{(\tau, \var)}_{k})_{mt}= 
\sqrt{\fr{t!}{m!}} f_k^{(\tau,\var}(m)\dl_{m,t+2k+2\var}
\eeq
for $k< 0$. 

\vspace{5mm}
\noindent
Finally we like to give two simple specific examples: Let $g=h(\alpha)$, 
$\var, k=0$ in (\r{54}) that is $\tau $ is positive integer. 
Taking $s, l=0$ in (\r{5a})  and  (\r{5b}) we get  
\beq\lb{lg}
P_\tau (\cosh \alpha )= \fr{1}{\sqrt{\pi}\cosh\fr{\alpha}{2}}
\sum_{n=0}^\infty \fr{\Gamma(n+\fr{1}{2})}{n!}\tanh^{2n}\fr{\alpha}{2}
F(-2n, 1+\tau; 1+n; -\sinh^2\fr{\alpha}{2})
\eeq
and
\beq
1= \sum_{n=0}^\tau   \fr{(-)^n!}{(n!)^2} \fr{(\tau +n)!}{(\tau- n)!}
\tanh^{2n}\fr{\alpha}{2}
F(-\tau, 1+\tau; 1+n; -\sinh^2\fr{\alpha}{2})
\eeq
respectively.  $P_\tau$ in (\r{lg}) is the Legendre function.

\end{document}